# Stochastic Comparisons of Series and Parallel Systems with Topp-Leone Generated Family of Distributions

**Ruby Chanchal**[1], **Vaishali Gupta**[2], **Amit Kumar Misra**[3]
Department of Statistics
School of Physical & Decision Sciences
Babasaheb Bhimrao Ambedkar University, Lucknow - 226 025

**Abstract**

In this article, we stochastically compare the series and parallel systems having Topp-Leone generated family of distributions. We consider that the lifetimes of the components of the systems have either the different shape parameters when the scale parameters are fixed or the different scale parameters when the shape parameters are fixed and established some ordering results with the help of vector majorization technique.

# 1 Introduction

Order statistics play a prominent role in statistics, applied probability, actuarial science, reliability theory and many other related fields. Let $X_{1:n} \leq X_{2:n} \leq \cdots \leq X_{n:n}$ represent the order statistics arising from the random variables $X_1, X_2, \ldots, X_n$. Also, let $X_{k:n}$ represents the $k^{\text{th}}$ order statistic which is related to the lifetimes of $(n-k+1)$-out-of-$n$ system. In particular, $X_{n:n}$ (when $k = n$) and $X_{1:n}$ (when $k = 1$) denote the lifetimes of parallel and series systems, respectively. It has been studied extensively in the literature and found to be fruitful in comparisons of lifetimes of series and parallel systems where the components' lifetimes are independent and identically distributed (i.i.d.). But in the case when the lifetimes of the components are non-i.i.d., the distribution theory becomes quite complicated. Because of this reason, hardly few results are available in the literature (see, for example, Balakrishnan and Rao (1998) and David and Nagaraja (2003), and references cited therein).

Many researchers have worked upon the stochastic comparisons between the lifetimes of different systems where the random lifetimes of components follow various lifetime distributions, for example, Dykstra et al. (1997) and Khaledi and Kochar (2000) studied the comparisons of parallel systems of heterogeneous exponential components. Khaledi and Kochar (2006), Fang and Tang (2014), and Torrado and Kochar (2015) considered the case of heterogeneous Weibull distributions. Moreover, several results have been derived for the heterogeneous generalized exponential distributions, gamma and Pareto type distributions (see, for examples, Balakrishnan and Zhao (2013), Balakrishnan et al. (2014), and Patra et al. (2018)). For recent references on the results related to the stochastic comparisons of this type, we refer to Barmalzan et al. (2016), Fang et al. (2016), Fang and Wang (2017), and Nadarajah et al. (2017). Apart from these types of comparisons, the families of distributions have also been considered. Some well known families of lifetime distributions are exponentiated Weibull (Mudholkar and Srivastava (1993)) and generalized exponential (Gupta and Kundu (1999)), etc. Recently, Kayal (2018) studied the stochastic comparisons of series and parallel systems with Kumaraswamy generalized family of distributions.

[1]E-mail: rubychanchal21@gmail.com
[2]E-mail: vaishali.gupta3091@gmail.com; Corresponding author
[3]E-mail: mishraamit31@gmail.com

It is important to mention that the notion of majorization is one of the useful tool to compare lifetimes of series and parallel systems. The concept of majorization deals with the diversity of components of vectors in $\mathbb{R}^n$. Let $\boldsymbol{x} = (x_1, x_2, \ldots, x_n)$ and $\boldsymbol{y} = (y_1, y_2, \ldots, y_n)$ be two real vectors. Further, let $x_{(1)} \leq x_{(2)} \leq \cdots \leq x_{(n)}$ and $y_{(1)} \leq y_{(2)} \leq \cdots \leq y_{(n)}$ denote the increasing order of the components of $\boldsymbol{x}$ and $\boldsymbol{y}$, respectively. The vector $\boldsymbol{x}$ is said to be majorized by vector $\boldsymbol{y}$ (written as $\boldsymbol{x} \overset{\text{m}}{\preceq} \boldsymbol{y}$) if $\sum_{k=1}^{i} x_{(k)} \geq \sum_{k=1}^{i} y_{(k)}$, $i = 1, \ldots, n-1$, and $\sum_{k=1}^{n} x_{(k)} = \sum_{k=1}^{n} y_{(k)}$. Thus, $\boldsymbol{x} \overset{\text{m}}{\preceq} \boldsymbol{y}$ tells that the components of $\boldsymbol{x}$ are less dispersed as compared to those of the $\boldsymbol{y}$. Majorization has played an important role to study several stochastic orders in various fields such as reliability theory, economics, mathematics, statistics and so on. For more details of majorization and its applications, one may refer to Marshall et al. (2011).

In this paper, we consider the stochastic comparisons of series and parallel systems with respect to the likelihood ratio order, the hazard rate order, and the usual stochastic order using vector majorization technique, where the components of the systems follow Topp-Leone generated family of distributions. This distribution was recently given by Rezaei et al. (2017) as a generalization of Topp and Leone's distribution (see, Topp and Leone (1955)). It has the property to model bathtub shaped hazard rates depending upon the values of parameters and it can be used for lifetime modelling. For more applications of this distribution, one may refer to Rezaei et al. (2017). Let $X$ be a random variable following the Topp-Leone generated (TL-$G$) family of distribution. The probability density function (p.d.f.) and the cumulative distribution function (c.d.f.) of $X$ are given by

$$f(x; \alpha, \theta, \xi) = 2\alpha\theta g(x; \xi) G(x; \xi)^{\theta\alpha - 1} (1 - G(x; \xi)^{\theta})(2 - G(x; \xi)^{\theta})^{\alpha - 1}, \quad x \geq 0,\ \theta, \alpha > 0 \tag{1.1}$$

and

$$F(x; \alpha, \theta, \xi) = (G(x; \xi)^{\theta}(2 - G(x; \xi)^{\theta}))^{\alpha}, \quad x \geq 0,\ \theta, \alpha > 0, \tag{1.2}$$

respectively, where $\theta$ is the scale parameter, $\alpha$ is the shape parameter, $G(x; \xi)$ is a baseline c.d.f., $g(x; \xi)$ is the corresponding p.d.f., and $\xi$ contains the parameters which specify the baseline distribution. For convenience, we use the notion $X \sim$TL-$G(\alpha, \theta, \xi)$.

Note that, Topp-Leone's distribution is a particular case of TL-$G$ family of distribution, when $G(x;\xi)$ be the $U(0,1)$ and $\theta = 1$. Various cases of TL-$G(\alpha,\theta,\xi)$ can be obtained by using $G(x;\xi)$ as a parent distribution function such as exponential, normal, log-logistic, and gamma, and we get TL-exponential, TL-normal, TL-log-logistic, and TL-gamma distributions, respectively.

The outline of the paper is as follows. In Section 2, we briefly present some definitions of stochastic orders and majorization, and discuss some lemmas which are useful to derive our results. In Section 3, we provide some ordering results with respect to the usual stochastic and the likelihood ratio orderings for the comparisons of parallel systems and with respect to the usual stochastic and the hazard rate orderings for the comparisons of series systems with TL-$G$ distributed components.

Throughout the paper, we use the notion $\mathbb{R} = (-\infty,\infty)$, $\mathbb{R}^+ = [0,\infty)$, and $\mathbb{R}^n = (-\infty,\infty)^n$. The terms increasing and decreasing are used for non-decreasing and non-increasing, respectively.

# 2 Preliminaries

In this section, we first recall some definitions of stochastic orders and vector majorization. Then, we discuss some useful lemmas which are essential to develop the results conferred in the next section.

## 2.1 Stochastic Orders

Let $X$ and $Y$ be two random variables having the common support $\mathbb{R}^+$ with the distribution functions $F_X(\cdot)$ and $F_Y(\cdot)$, the probability density functions $f_X(\cdot)$ and $f_Y(\cdot)$, and the hazard rate functions $r_X(\cdot)$ and $r_Y(\cdot)$, respectively. Let $\bar{F}_X(\cdot) = 1 - F_X(\cdot)$ and $\bar{F}_Y(\cdot) = 1 - F_Y(\cdot)$ be the survival functions of $X$ and $Y$, respectively.

**Definition 2.1.** *$X$ is said to be smaller than $Y$ in the*

*(i) usual stochastic order (written as $X \leq_{st} Y$) if $F_Y(x) \leq F_X(x)$, $\forall x \in \mathbb{R}^+$;*

*(ii) hazard rate order (written as $X \leq_{hr} Y$) if $\bar{F}_Y(x)/\bar{F}_X(x)$ is increasing in $x \in \mathbb{R}^+$, or equivalently if $r_X(x) \geq r_Y(x)$, $\forall x \in \mathbb{R}^+$;*

*(iii) likelihood ratio order (written as $X \leq_{lr} Y$) if $f_Y(x)/f_X(x)$ is increasing in $x \in \mathbb{R}^+$.*

The following implications are well-known.

$$X \leq_{\text{lr}} Y \Rightarrow X \leq_{\text{hr}} Y \Rightarrow X \leq_{\text{st}} Y.$$

For an extensive details on various stochastic orders, see, Shaked and Shanthikumar (2007).

## 2.2 Majorization

In this subsection, we briefly present some definitions and results related to the vector majorization which are very useful in dealing with various inequalities while comparing the order statistics. Let $\boldsymbol{x} = (x_1, x_2, \ldots, x_n)$ and $\boldsymbol{y} = (y_1, y_2, \ldots, y_n)$ be two real vectors from $\mathbb{R}^n$. Further, let $x_{(1)} \leq x_{(2)} \leq \cdots \leq x_{(n)}$ and $y_{(1)} \leq y_{(2)} \leq \cdots \leq y_{(n)}$ denote the increasing order of the components of $\boldsymbol{x}$ and $\boldsymbol{y}$, respectively.

**Definition 2.2.** *$\boldsymbol{x}$ is said to be*

*(i) majorized by $\boldsymbol{y}$ (written as $\boldsymbol{x} \overset{m}{\preceq} \boldsymbol{y}$) if*

$$\sum_{k=1}^{i} x_{(k)} \geq \sum_{k=1}^{i} y_{(k)},\ i = 1, \ldots, n-1,\ \text{and} \sum_{k=1}^{n} x_{(k)} = \sum_{k=1}^{n} y_{(k)}.$$

*(ii) weakly submajorized by $\boldsymbol{y}$ (written as $\boldsymbol{x} \preceq_w \boldsymbol{y}$) if*

$$\sum_{k=i}^{n} x_{(k)} \leq \sum_{k=i}^{n} y_{(k)},\ i = 1, \ldots, n.$$

From the Definition 2.2, it is easy to verify that $\boldsymbol{x} \overset{\mathrm{m}}{\preceq} \boldsymbol{y}$ implies $\boldsymbol{x} \preceq_{\mathrm{w}} \boldsymbol{y}$. For a comprehensive study on majorization, we refer the reader to Marshall et al. (2011) .

**Definition 2.3.** *A real valued function $\psi$ defined on a set $\mathbb{A} \subset \mathbb{R}^n$ is said to be Schur-convex (Schur-concave) on $\mathbb{A}$ if*

$$\boldsymbol{x} \overset{m}{\preceq} \boldsymbol{y} \implies \psi(\boldsymbol{x}) \leq (\geq)\, \psi(\boldsymbol{y}) \quad \text{for } \boldsymbol{x}, \boldsymbol{y} \in \mathbb{A}.$$

**Lemma 2.1.** *(Marshall et al., 2011, p.84) Let $I \subset \mathbb{R}$ be an open interval and let $\psi : I^n \to \mathbb{R}$ be continuously differentiable. Necessary and sufficient conditions for $\psi(\cdot)$ to be Schur-concave on $I^n$ are: $\psi(\cdot)$ is symmetric on $I^n$, and*

$$(x_k - x_l)\left(\frac{\partial \psi(\boldsymbol{x})}{\partial x_k} - \frac{\partial \psi(\boldsymbol{x})}{\partial x_l}\right) \leq 0,\ \text{for all } k \neq l.$$

**Lemma 2.2.** *(Marshall et al., 2011, p.87) A real valued function $\psi$ defined on $\mathbb{A} \subset \mathbb{R}^n$ satisfies $\boldsymbol{x} \preceq_w \boldsymbol{y}$ on $\mathbb{A}$ implies $\psi(\boldsymbol{x}) \leq \psi(\boldsymbol{y})$ if, and only if, $\psi$ is increasing and Schur-convex on $\mathbb{A}$.*

**Lemma 2.3.** *(Balakrishnan et al., 2014) Let a function $\tau : (0, \infty) \to (0, \infty)$ be defined as*

$$\tau(\alpha) = \frac{\alpha t^{\alpha-1}}{1 - t^{\alpha}}.$$

*Then, $\tau(\alpha)$ is convex in $\alpha$, for any $0 < t < 1$.*

**Lemma 2.4.** *(Marshall et al., 2011, Proposition C.1., p.92) If $I \subset \mathbb{R}$ is an open interval and $h : I \to \mathbb{R}$ is convex, then $\psi(\boldsymbol{x}) = \sum_{i=1}^{n} h(x_i)$ is Schur-convex on $I^n$. Consequently, $\boldsymbol{x} \overset{m}{\preceq} \boldsymbol{y}$ on $I^n$ implies $\psi(\boldsymbol{x}) \leq \psi(\boldsymbol{y})$.*

# 3 Main Results and Conclusions

In this section, we compare the lifetimes of series and parallel systems having independent TL-$G$ distributed components. These results are presented with heterogeneity in one parameter while another is fixed.

The following theorem deals with the hazard rate ordering of series systems when the parameter $\boldsymbol{\alpha} = (\alpha_1, \alpha_2, \ldots, \alpha_n)$ varies.

**Theorem 3.1.** *Let $X_1, X_2, \ldots, X_n$ and $Y_1, Y_2, \ldots, Y_n$ be the two sets of independent random variables with $X_k \sim$ TL-$G(\alpha_k, \theta, \xi)$ and $Y_k \sim$ TL-$G(\alpha_k^*, \theta, \xi)$ for $k = 1, 2, \ldots, n$, respectively. Then, for fixed $\theta > 0$ and for any fixed $\xi$, we have*

$$\boldsymbol{\alpha}^* = (\alpha_1^*, \alpha_2^*, \ldots, \alpha_n^*) \overset{m}{\preceq} (\alpha_1, \alpha_2, \ldots, \alpha_n) = \boldsymbol{\alpha} \quad \Longrightarrow \quad X_{1:n} \leq_{hr} Y_{1:n}.$$

*Proof.* It is well defined for a series system that the sum of the hazard rate functions of each components is equal to the hazard rate function of the system. Therefore, for $x > 0$, the hazard rate function of $X_{1:n}$ is given by

$$\begin{aligned}
r_{X_{1:n}}(x) &= \sum_{k=1}^{n} \frac{f(x;\alpha_k,\theta,\xi)}{\bar{F}(x;\alpha_k,\theta,\xi)} \\
&= \sum_{k=1}^{n} \frac{2\alpha_k\theta g(x;\xi)G(x;\xi)^{\theta\alpha_k-1}(1-G(x;\xi)^{\theta})(2-G(x;\xi)^{\theta})^{\alpha_k-1}}{1-G(x;\xi)^{\theta\alpha_k}(2-G(x;\xi)^{\theta})^{\alpha_k}} \\
&= 2\theta g(x;\xi)G(x;\xi)^{\theta-1}(1-G(x;\xi)^{\theta}) \sum_{k=1}^{n} \frac{\alpha_k G(x;\xi)^{\theta(\alpha_k-1)}(2-G(x;\xi)^{\theta})^{\alpha_k-1}}{1-G(x;\xi)^{\theta\alpha_k}(2-G(x;\xi)^{\theta})^{\alpha_k}} \\
&= 2\theta g(x;\xi)G(x;\xi)^{\theta-1}(1-G(x;\xi)^{\theta}) \sum_{k=1}^{n} \frac{\alpha_k (G(x;\xi)^{\theta}(2-G(x;\xi)^{\theta}))^{\alpha_k-1}}{1-(G(x;\xi)^{\theta}(2-G(x;\xi)^{\theta}))^{\alpha_k}} \\
&= 2\theta g(x;\xi)G(x;\xi)^{\theta-1}(1-G(x;\xi)^{\theta}) \sum_{k=1}^{n} z(\alpha_k),
\end{aligned}$$

where, for fixed $x > 0$, $\theta > 0$, and for any fixed $\xi$,

$$z(\alpha) = \frac{\alpha (G(x;\xi)^{\theta}(2-G(x;\xi)^{\theta}))^{\alpha-1}}{1-(G(x;\xi)^{\theta}(2-G(x;\xi)^{\theta}))^{\alpha}}, \quad \alpha > 0.$$

On taking $t = G(x;\xi)^{\theta}(2 - G(x;\xi)^{\theta})$ and using Lemma 2.3, it follows that $z(\alpha)$ is convex in $\alpha$. Now, on using Lemma 2.4, we conclude that $\sum_{k=1}^{n} z(\alpha_k)$ is Schur-convex on $(0,\infty)^n$, which implies that if $\boldsymbol{\alpha}^* \overset{m}{\preceq} \boldsymbol{\alpha}$, then $r_{X_{1:n}}(x) \geq r_{Y_{1:n}}(x)$. Hence the theorem follows. □

Now, we present the following example to discuss the above theorem.

**Example 3.1.** Consider that $G(x;\xi) = 1 - e^{-x},\ x \geq 0$. Let $X_1, X_2$ and $Y_1, Y_2$ be the two sets of independent random variables with $X_k \sim$ TL-G$(\alpha_k, \theta, \xi)$ and $Y_k \sim$ TL-G$(\alpha_k^*, \theta, \xi)$ for $k = 1, 2$, respectively. Assume $\alpha_1 = 1, \alpha_2 = 9,\ \alpha_1^* = 4, \alpha_2^* = 6$, and $\theta = 0.5$. Clearly, $(\alpha_1^*, \alpha_2^*) \overset{m}{\preceq} (\alpha_1, \alpha_2)$, and therefore, using Theorem 3.1, we have $X_{1:2} \leq_{\text{hr}} Y_{1:2}$. This can also be concluded from the Figure 1(a) where we plot $r_{X_{1:2}}(x) - r_{Y_{1:2}}(x)$ which is non-negative for $x \geq 0$. □

One may be interested to know whether Theorem 3.1 can be extended to the likelihood ratio order. The following counterexample shows that the result will not hold.

**Counterexample 3.1.** Continuing with the Example 3.1, if we plot $\frac{f_{Y_{1:2}}(x)}{f_{X_{1:2}}(x)}$, we get the Figure 1(b), which shows that as the value of $x$ increases, the ratio first increases and then decreases. Hence the result in Theorem 3.1 cannot be extended to the likelihood ratio order. □

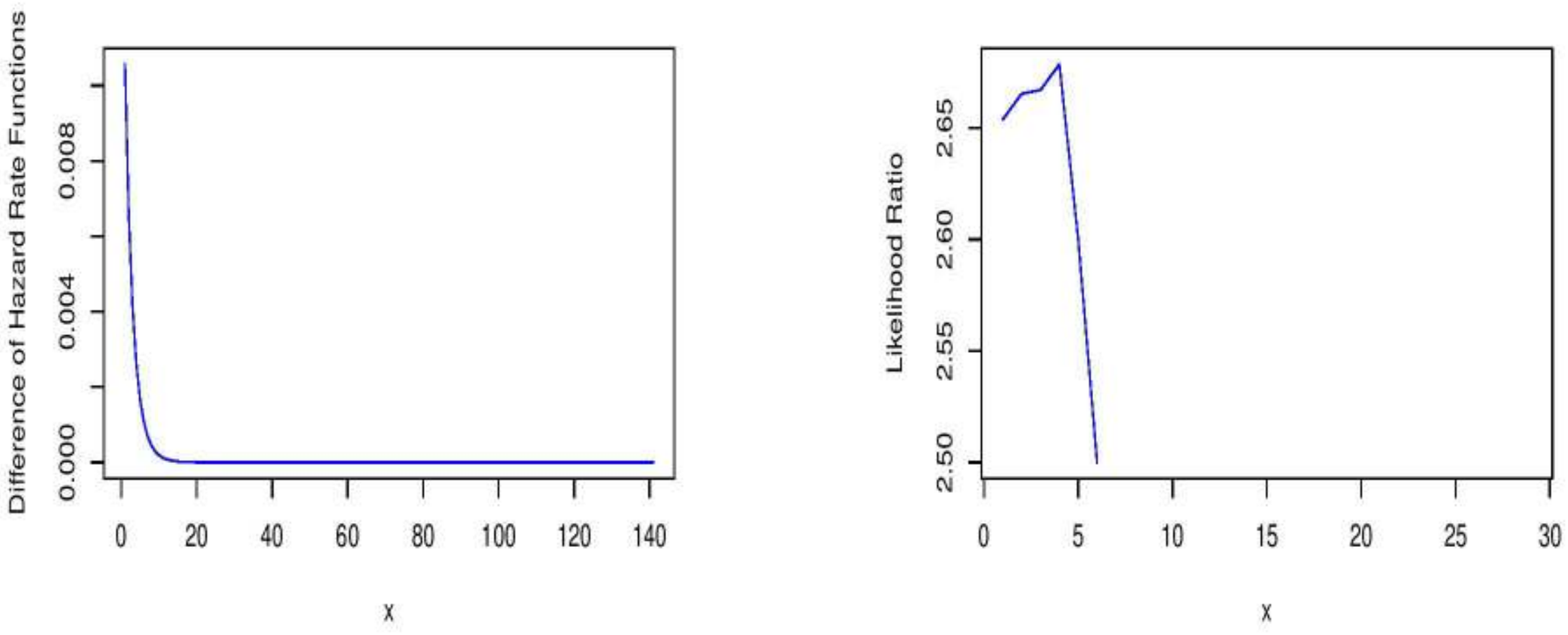

(a) Plot of $r_{X_{1:2}}(x) - r_{Y_{1:2}}(x)$

(b) Plot of $\frac{f_{Y_{1:2}}(x)}{f_{X_{1:2}}(x)}$

Figure 1

In the next theorems, we provide the stochastic comparisons of parallel systems in the sense of the usual stochastic order and the likelihood ratio order when the parameters $\boldsymbol{\theta} = (\theta_1, \theta_2, \ldots, \theta_n)$ and $\boldsymbol{\alpha} = (\alpha_1, \alpha_2, \ldots, \alpha_n)$ vary.

**Theorem 3.2.** *Let $X_1, X_2, \dots, X_n$ and $Y_1, Y_2, \dots, Y_n$ be the two sets of independent random variables with $X_k \sim TL\text{-}G(\alpha, \theta_k, \xi)$ and $Y_k \sim TL\text{-}G(\alpha, \theta_k^*, \xi)$ for $k = 1, 2, \dots, n$, respectively. Then, for fixed $\alpha > 0$ and for any fixed $\xi$, we have*

$$\boldsymbol{\theta} = (\theta_1, \theta_2, \dots, \theta_n) \preceq_w (\theta_1^*, \theta_2^*, \dots, \theta_n^*) = \boldsymbol{\theta}^* \quad \Longrightarrow \quad X_{n:n} \leq_{st} Y_{n:n}.$$

*Proof.* The distribution function of $X_{n:n}$ is given by

$$F_{X_{n:n}}(x) = \prod_{i=1}^{n} G(x;\xi)^{\theta_i \alpha}(2 - G(x;\xi)^{\theta_i})^{\alpha} = \varphi(\boldsymbol{\theta}), \text{ say.} \tag{3.1}$$

Clearly, $\varphi : (0, \infty)^n \to \mathbb{R}$ is a symmetric function on $(0, \infty)^n$. Partially differentiating $\varphi(\boldsymbol{\theta})$ with respect to $\theta_k$, we get

$$\begin{aligned}
\frac{\partial \varphi(\boldsymbol{\theta})}{\partial \theta_k} &= \left[\prod_{\substack{i=1\\ i\neq k}}^{n} G(x;\xi)^{\theta_i \alpha}(2 - G(x;\xi)^{\theta_i})^{\alpha}\right]\Big( - \alpha G(x;\xi)^{\theta_k(\alpha+1)}(2 - G(x;\xi)^{\theta_k})^{\alpha-1} \ln G(x;\xi) \\
&\quad + (2 - G(x;\xi)^{\theta_k})^{\alpha} G(x;\xi)^{\theta_k \alpha} \ln G(x;\xi)^{\alpha}\Big) \\
&= \left[\prod_{\substack{i=1\\ i\neq k}}^{n} G(x;\xi)^{\theta_i \alpha}(2 - G(x;\xi)^{\theta_i})^{\alpha}\right] \alpha G(x;\xi)^{\theta_k \alpha}(2 - G(x;\xi)^{\theta_k})^{\alpha} \ln G(x;\xi) \\
&\quad \times \left(1 - \frac{G(x;\xi)^{\theta_k}}{2 - G(x;\xi)^{\theta_k}}\right) \\
&= 2\alpha F_{X_{n:n}}(x) \ln G(x;\xi) \left(\frac{1 - G(x;\xi)^{\theta_k}}{2 - G(x;\xi)^{\theta_k}}\right).
\end{aligned}$$

It is easy to see that $\frac{\partial \varphi(\boldsymbol{\theta})}{\partial \theta_k} < 0$. Therefore, $\varphi(\boldsymbol{\theta})$ is decreasing in $\theta_k$. For $\theta_k \neq \theta_l$, we have

$$\begin{aligned}
(\theta_k - \theta_l)\left(\frac{\partial \varphi(\boldsymbol{\theta})}{\partial \theta_k} - \frac{\partial \varphi(\boldsymbol{\theta})}{\partial \theta_l}\right) &= 2\alpha(\theta_k - \theta_l) F_{X_{n:n}}(x) \ln G(x;\xi) \\
&\quad \times \left(\frac{1 - G(x;\xi)^{\theta_k}}{2 - G(x;\xi)^{\theta_k}} - \frac{1 - G(x;\xi)^{\theta_l}}{2 - G(x;\xi)^{\theta_l}}\right) \\
&= 2\alpha(\theta_k - \theta_l) F_{X_{n:n}}(x) \ln G(x;\xi) \\
&\quad \times \left(\frac{G(x;\xi)^{\theta_l} - G(x;\xi)^{\theta_k}}{(2 - G(x;\xi)^{\theta_k})\,(2 - G(x;\xi)^{\theta_l})}\right) \\
&\leq 0.
\end{aligned}$$

On using Lemma 2.1, $\varphi(\boldsymbol{\theta})$ is Schur-concave in $\boldsymbol{\theta}$. Thus, $-\varphi(\boldsymbol{\theta})$ is increasing in $\theta_k$ and Schur-convex in $\boldsymbol{\theta}$. Now, using Lemma 2.2, it follows that $\boldsymbol{\theta} \preceq_w \boldsymbol{\theta}^*$ implies $-\varphi(\boldsymbol{\theta}) \leq -\varphi(\boldsymbol{\theta}^*)$, or equivalently, $\varphi(\boldsymbol{\theta}^*) \leq \varphi(\boldsymbol{\theta})$. Therefore, $F_{Y_{n:n}}(x) \leq F_{X_{n:n}}(x)$ and hence $X_{n:n} \leq_{\text{st}} Y_{n:n}$. □

The following corollary is an immediate consequence of Theorem 3.2.

**Corollary 3.1.** *Let $X_1, X_2, \ldots, X_n$ and $Y_1, Y_2, \ldots, Y_n$ be the two sets of independent random variables with $X_k \sim TL\text{-}G(\alpha, \theta_k, \xi)$ and $Y_k \sim TL\text{-}G(\alpha, \theta_k^*, \xi)$, $k = 1, 2, \ldots, n$, respectively. Then, for fixed $\alpha > 0$ and for any fixed $\xi$, we have*

$$\boldsymbol{\theta} \overset{m}{\preceq} \boldsymbol{\theta}^* \implies X_{n:n} \leq_{st} Y_{n:n}.$$

The following example illustrates the result established in Theorem 3.2.

**Example 3.2.** Consider that $G(x; \xi) = 1 - e^{-x}$, $x \geq 0$. Let $X_1, X_2$ and $Y_1, Y_2$ be the two sets of independent random variables with $X_k \sim$ TL-G$(\alpha, \theta_k, \xi)$ and $Y_k \sim$ TL-G$(\alpha, \theta_k^*, \xi)$ for $k = 1, 2$, respectively. Assume $\theta_1 = 0.1$, $\theta_2 = 0.4$, $\theta_1^* = 0.2$, $\theta_2^* = 0.5$, and $\alpha = 0.5$. Clearly, $(\theta_1, \theta_2) \preceq_{\text{w}} (\theta_1^*, \theta_2^*)$, and therefore, using Theorem 3.2, we have $X_{2:2} \leq_{\text{st}} Y_{2:2}$. This can also be seen from the Figure 2(a) where we plot $F_{X_{2:2}}(x) - F_{Y_{2:2}}(x)$ which is non-negative for $x \geq 0$. □

The following counterexample shows that the result in Theorem 3.2 may not hold for the likelihood ratio order.

**Counterexample 3.2.** Continuing with the Example 3.2, if we plot $\frac{f_{Y_{2:2}}(x)}{f_{X_{2:2}}(x)}$, we get the Figure 2(b), which shows that as the value of $x$ increases, the ratio first increases and then decreases. Thus, the result in Theorem 3.2 cannot be extended to the likelihood ratio order. □

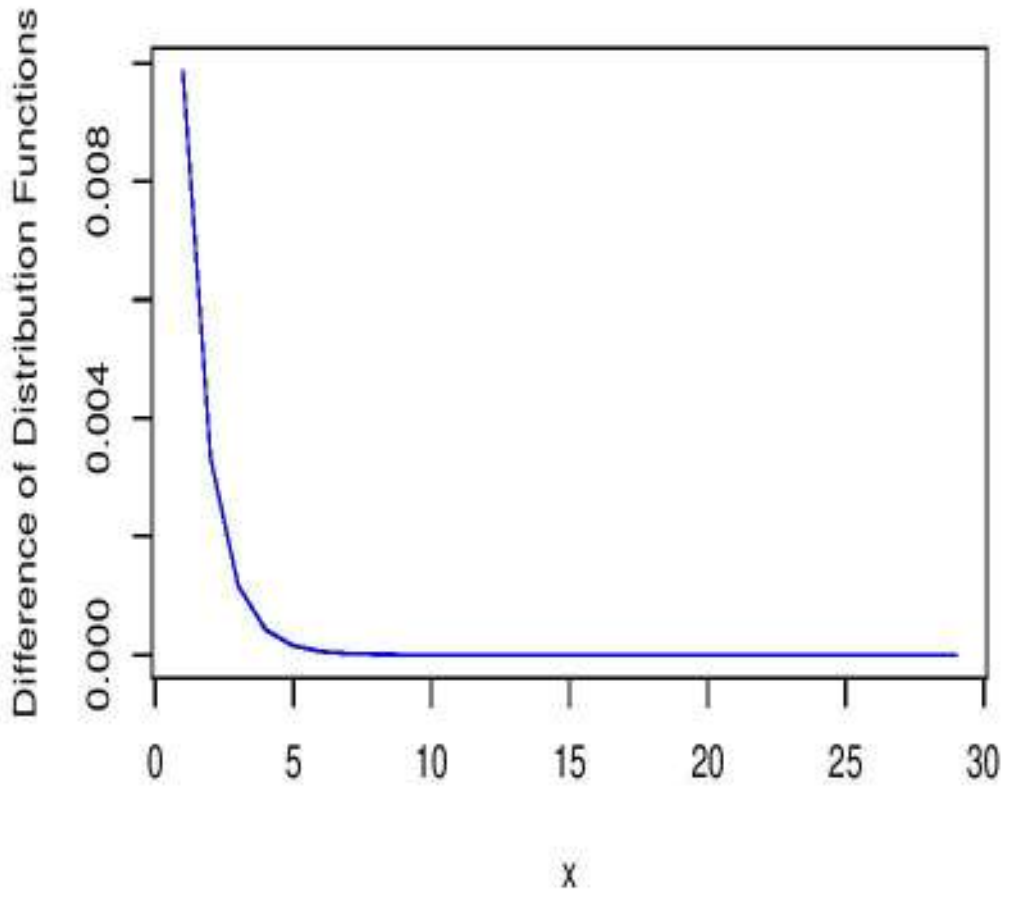

(a) Plot of $F_{X_{2:2}}(x) - F_{Y_{2:2}}(x)$

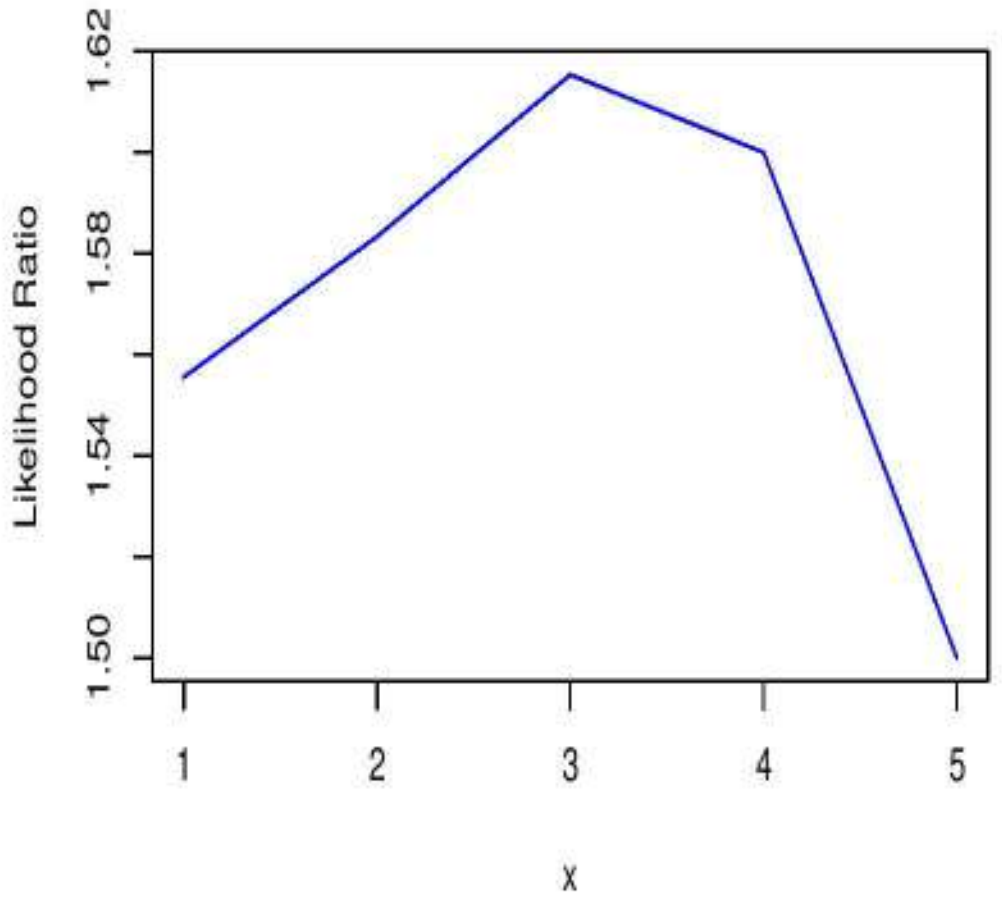

(b) Plot of $\frac{f_{Y_{2:2}}(x)}{f_{X_{2:2}}(x)}$

Figure 2

Next result is a generalization of Theorem 3.2 to a wide range of scale parameters.

**Theorem 3.3.** *Let $X_1, X_2, \ldots, X_n$ and $Y_1, Y_2, \ldots, Y_n$ be the two sets of independent random variables with $X_k \sim TL\text{-}G(\alpha, \theta_k, \xi)$ and $Y_k \sim TL\text{-}G(\alpha, \theta_k^*, \xi)$, $k = 1, 2, \ldots, n$, respectively. For fixed $\alpha > 0$ and for any fixed $\xi$, if $(\theta_1, \theta_2, \ldots, \theta_n) \leq (\theta_1^*, \theta_2^*, \ldots, \theta_n^*)$, that is, $\theta_k \leq \theta_k^*$, $k = 1, 2, \ldots, n$, then $X_{n:n} \leq_{st} Y_{n:n}$.*

*Proof.* In the proof of Theorm 3.2, we have shown that $\varphi(\boldsymbol{\theta})$, given by (3.1), is decreasing in each $\theta_k$, $k \in \{1, 2, \ldots, n\}$. Therefore, $\varphi(\boldsymbol{\theta}^*) \leq \varphi(\boldsymbol{\theta})$, or equivalently, $F_{Y_{n:n}}(x) \leq F_{X_{n:n}}(x)$. Hence the required result follows. □

**Theorem 3.4.** *Let $X_1, X_2, \ldots, X_n$ and $Y_1, Y_2, \ldots, Y_n$ be the two sets of independent random variables with $X_k \sim TL\text{-}G(\alpha_k, \theta, \xi)$ and $Y_k \sim TL\text{-}G(\alpha_k^*, \theta, \xi)$, $k = 1, 2, \ldots, n$, respectively. Then, for fixed $\theta > 0$ and for any fixed $\xi$, $X_{n:n} \leq_{lr} Y_{n:n}$ if, and only if, $\sum_{k=1}^n \alpha_k \leq \sum_{k=1}^n \alpha_k^*$.*

*Proof.* It is easy to verify that the distribution function and the probability density function of $X_{n:n}$ are given by

$$F_{X_{n:n}}(x) = G(x;\xi)^{\theta \sum_{k=1}^n \alpha_k} (2 - G(x;\xi)^\theta)^{\sum_{k=1}^n \alpha_k}, \ x \geq 0,$$

and

$$f_{X_{n:n}}(x) = 2\theta \left( \sum_{k=1}^n \alpha_k \right) g(x;\xi) G(x;\xi)^{\theta \sum_{k=1}^n \alpha_k - 1} (2 - G(x;\xi)^\theta)^{\sum_{k=1}^n \alpha_k} \left( \frac{1 - G(x;\xi)^\theta}{2 - G(x;\xi)^\theta} \right), \ x > 0,$$

respectively. Similarly, the probability density function of $Y_{n:n}$ is given by

$$f_{Y_{n:n}}(x) = 2\theta \left( \sum_{k=1}^n \alpha_k^* \right) g(x;\xi) G(x;\xi)^{\theta \sum_{k=1}^n \alpha_k^* - 1} (2 - G(x;\xi)^\theta)^{\sum_{k=1}^n \alpha_k^*} \left( \frac{1 - G(x;\xi)^\theta}{2 - G(x;\xi)^\theta} \right), \ x > 0.$$

Then, for $x > 0$,

$$\frac{f_{Y_{n:n}}(x)}{f_{X_{n:n}}(x)} = \left( \frac{\sum_{k=1}^n \alpha_k^*}{\sum_{k=1}^n \alpha_k} \right) (G(x;\xi)^\theta (2 - G(x;\xi)^\theta))^{\sum_{k=1}^n \alpha_k^* - \sum_{k=1}^n \alpha_k}.$$

Using the fact that $G(x;\xi)^\theta$ is an increasing function of $x$ and the observation that $y(2-y)$ is an increasing function of $y \in (0,1)$, we conclude that the ratio $\frac{f_{Y_{n:n}}(x)}{f_{X_{n:n}}(x)}$ is increasing in $x$ if, and only if, $\sum_{k=1}^n \alpha_k \leq \sum_{k=1}^n \alpha_k^*$, which proves the desired result. □

Till now we have derived the results when TL-$G$ family of distributions have the same baseline distributions. Now, we present the results for the case when TL-$G$ family of distributions have different baseline distributions. Let $X_1^*$ and $X_2^*$ be the two random variables having the cumulative distribution functions $G_1(\cdot)$ and $G_2(\cdot)$, respectively. Also, assume that $U_1, U_2, \ldots, U_n$ and $V_1, V_2, \ldots, V_n$ be the two sets of independent random variables following the TL-$G$ family of distributions with baseline distributions $G_1(\cdot)$ and $G_2(\cdot)$, respectively, and we denote $U_k \sim$ TL-G$(\alpha_k, \theta, G_1)$ and $V_k \sim$ TL-G$(\alpha_k^*, \theta, G_2)$ for $k = 1, \dots, n$. The survival functions of $U_{1:n}$ and $V_{1:n}$ are respectively given by

$$\bar{F}_{U_{1:n}}(x) = \prod_{k=1}^{n} \left[1 - \left((G_1(x))^\theta \left(2 - (G_1(x))^\theta\right)\right)^{\alpha_k}\right], \ x \geq 0$$

and

$$\bar{F}_{V_{1:n}}(x) = \prod_{k=1}^{n} \left[1 - \left((G_2(x))^\theta \left(2 - (G_2(x))^\theta\right)\right)^{\alpha_k^*}\right], \ x \geq 0.$$

The following theorem provides the conditions under which $U_{1:n} \leq_{\text{st}} V_{1:n}$.

**Theorem 3.5.** *Let $U_1, U_2, \dots, U_n$ and $V_1, V_2, \dots, V_n$ be the two set of independent random variables with $U_k \sim$ TL-G$(\alpha_k, \theta, G_1)$ and $V_k \sim$ TL-G$(\alpha_k^*, \theta, G_2)$ for $k = 1, 2, \dots, n$, respectively, and let $\boldsymbol{\alpha}^* \overset{m}{\preceq} \boldsymbol{\alpha}$. Then, for fixed $\theta > 0$,*

$$X_1^* \leq_{st} X_2^* \implies U_{1:n} \leq_{st} V_{1:n}.$$

*Proof.* Let $Z_1, Z_2, \dots, Z_n$ be the set of independent random variable with $Z_k \sim$ TL-G$(\alpha_k^*, \theta, G_1)$ for $k = 1, 2, \dots, n$. On using Theorem 3.1, we have $U_{1:n} \leq_{\text{hr}} Z_{1:n}$, which implies that $U_{1:n} \leq_{\text{st}} Z_{1:n}$. Also, the survival function of $Z_{1:n}$ is given by

$$\bar{F}_{Z_{1:n}}(x) = \prod_{k=1}^{n} \left[1 - \left((G_1(x))^\theta \left(2 - (G_1(x))^\theta\right)\right)^{\alpha_k^*}\right], \ x \geq 0.$$

Since $X_1^* \leq_{\text{st}} X_2^*$ implies that $G_2(x) \leq G_1(x)$ for all $x \geq 0$, which further implies that $(G_2(x))^\theta \leq (G_1(x))^\theta$ for all $x \geq 0$. Now, using the observation that $y(2-y)$ is an increasing function of $y \in (0, 1)$, we have

$$\prod_{k=1}^{n} \left[1 - \left((G_1(x))^\theta \left(2 - (G_1(x))^\theta\right)\right)^{\alpha_k^*}\right] \leq \prod_{k=1}^{n} \left[1 - \left((G_2(x))^\theta \left(2 - (G_2(x))^\theta\right)\right)^{\alpha_k^*}\right], \forall x \geq 0, \tag{3.2}$$

i.e., $\bar{F}_{Z_{1:n}}(x) \leq \bar{F}_{V_{1:n}}(x)$ for all $x \geq 0$. Therefore, $Z_{1:n} \leq_{\text{st}} V_{1:n}$. Thus, we have $U_{1:n} \leq_{\text{st}} Z_{1:n} \leq_{\text{st}} V_{1:n}$. Hence the result follows. ☐

The following theorem provides the sufficient conditions for the comparision of parallel systems.

**Theorem 3.6.** *Let $W_1, W_2, \dots, W_n$ and $W_1^*, W_2^*, \dots, W_n^*$ be the two sets of independent random variables with $W_k \sim$ TL-G$(\alpha, \theta_k, G_1)$ and $W_k^* \sim$ TL-G$(\alpha, \theta_k^*, G_2)$ for $k = 1, 2, \dots, n$, respectively. For fixed $\alpha > 0$, if*

*(i) $\boldsymbol{\theta} \preceq_w \boldsymbol{\theta}^*$, then $X_1^* \leq_{st} X_2^*$ implies $W_{n:n} \leq_{st} W_{n:n}^*$;*

*(ii) $(\theta_1, \theta_2, \dots, \theta_n) \leq (\theta_1^*, \theta_2^*, \dots, \theta_n^*)$, i.e., if $\theta_k \leq \theta_k^*$, $k = 1, 2, \dots, n$, then $X_1^* \leq_{st} X_2^*$ implies $W_{n:n} \leq_{st} W_{n:n}^*$.*

*Proof.* Let $Z_1^*, Z_2^*, \ldots, Z_n^*$ be the set of independent random variable with $Z_k \sim \text{TL} - G(\alpha, \theta_k^*, G_1)$ for $k = 1, 2, \ldots, n$. On using the arguments similar to that used in the proof of Theorem 3.5, the part $(i)$ and part $(ii)$ follows from the Theorem 3.2 and Theorem 3.3, respectively. $\square$